\documentclass[final,3p,times,10pt]{elsarticle}
\usepackage{graphicx,amsmath,amssymb,txfonts}
\usepackage{algorithm}
\usepackage[noend]{algpseudocode}
\usepackage{algorithmicx,algorithm}
\usepackage{color}
\usepackage[colorlinks]{hyperref}
\usepackage{hypernat}
\usepackage{multirow}
\usepackage{enumerate}
\usepackage{booktabs}
\usepackage{epstopdf}
\usepackage{mathrsfs}%% \mathscr
\usepackage{subfigure}
\numberwithin{equation}{section}
\graphicspath{{Figs/}}

\biboptions{sort&compress}

\makeatletter

\newcommand{\Rmnum}[1]{\expandafter\@slowromancap\romannumeral #1@}
\makeatother

\newtheorem{remark}{\textbf{Remark}}

\newenvironment{proof}{\noindent{\bf Proof:}}{\hfill\fbox{}\vspace*{1mm}}

\usepackage{cases}
\usepackage{color}
\newcommand\bcase{\begin{numcases}{}}
\newcommand\ecase{\end{numcases}}

\begin{document}
\begin{frontmatter}



\title{The energy method for high-order invariants in shallow water wave equations}
\author[ZSTU]{Qifeng Zhang}
\author[ZSTU]{Tong Yan}
\author[NUPT]{Guang-hua Gao}
%\author[SEU]{Zhi-zhong Sun}
\cortext[cor]{E-mail address: zhangqifeng0504@gmail.com (Q. Zhang), tyan0320@mails.zstu.edu.cn (Tong Yan),
 gaogh@njupt.edu.cn (G. Gao)%,  \\zzsun@seu.edu.cn (Z. Sun)
}
\address[ZSTU]{Department of Mathematics, Zhejiang Sci-Tech University, Hangzhou, 310018, China}
\address[NUPT]{Department of Mathematics, Nanjing University of Posts and Telecommunications, Nanjing, 210096, China}
%\address[SEU]{School of Mathematics, Southeast University, Nanjing, 210096}

\begin{abstract}
Third order dispersive evolution equations are widely adopted to model one-dimensional long waves and have extensive applications in fluid mechanics, plasma physics and nonlinear optics. Among them are the KdV equation, the Camassa--Holm equation and the Degasperis--Procesi equation. They share many common features such as complete integrability, Lax pairs and bi-Hamiltonian structure.
In this paper we revisit high-order invariants for these three types of shallow water wave equations by the energy method in combination of a skew-adjoint operator $(1-\partial_{xx})^{-1}$. Several applications to seek high-order invariants of the Benjamin-Bona-Mahony equation, the regularized long wave equation and the Rosenau equation are also presented.
\end{abstract}

 \end{frontmatter}

\section{Introduction}
\setcounter{equation}{0}
A family of third order dispersive evolution equations of the form
\begin{align}
u_t - \alpha^2u_{xxt} + \gamma u_{xxx} + c_0 u_x = (c_1u^2 + c_2u_x^2 + c_3uu_{xx})_x,\quad x\in R,\;t>0   \label{eq1.1}
\end{align}
frequently appeared in the simulation of the shallow water waves, see e.g., \cite{ELY2006}, where $\alpha$, $\gamma$ and $c_i$ $(i=0,1,2,3)$ are real constants; $u$ denotes a horizontal velocity field with the independent spatial variable $x$ and temporal variable $t$.

A typical such equation \eqref{eq1.1} with $\alpha^2 = c_0 = c_2 = c_3 = 0$, $c_1 = 2$, $\gamma=-2$ is \emph{the KdV equation}
\begin{align}
u_t -4 uu_x - 2  u_{xxx}  = 0, \quad x\in R,\;t>0,  \label{eq1.2}
\end{align}
which describes the unidirectional propagation of waves at the free surface of shallow water under the influence of gravity.
The first four invariants of \eqref{eq1.2} are respectively as (see e.g., \cite{Tao2002}, although there is a
minor typo in the coefficient of the fourth invariant, it does not affect the reading of this classic review)
\begin{align*}
	&M_1 = \int_R u \mathrm{d}x,  \quad
	M_2 = \int_R u^2 \mathrm{d}x,   \quad M_3 = \int_R \Big(u_x^2 - \frac{2}{3}u^3\Big)\mathrm{d}x,\quad
    M_4 = \int_R \Big( u_{xx}^2-\frac{10}{3}uu_x^2 +\frac{5}{9}u^4\Big) \mathrm{d}x.
\end{align*}

Taking $\alpha^2 = c_3 = 1$, $\gamma = c_0 = 0$, $c_1 = -\frac{3}{2}$, $c_2 = \frac{1}{2}$,
we have another example called \emph{the Camassa--Holm equation} \cite{CH1993}
\begin{align}
u_t - u_{xxt} + 3uu_x = 2u_xu_{xx} + uu_{xxx},\quad x\in R,\;t>0,   \label{eq1.3}
\end{align}
which models the unidirectional propagation of shallow water waves over a flat bottom.
The first three invariants are listed as follows
\begin{align*}
	&E_1 = \int_R (u-u_{xx})\mathrm{d}x,  \quad
	E_2 = \frac{1}{2}\int_R (u^2 + u_x^2) \mathrm{d}x,   \quad
	E_3 = \frac{1}{2} \int_R u(u^2 + u_x^2) \mathrm{d}x.
\end{align*}

The third example by assigning $\alpha^2 = c_2 = c_3 = 1$, $\gamma = c_0 = 0$, $c_1 = -2$ is called \emph{the Degasperis--Procesi equation}
\begin{align}
	u_t - u_{xxt} + 4uu_x = 3u_xu_{xx} + uu_{xxx},\quad x\in R,\;t>0,  \label{eq1.4}
\end{align}
which can be regarded as a model for nonlinear shallow water dynamics \cite{DP1999}.
The frequently discussed invariants are
\begin{align*}
	&H_1 = \int_R (u-u_{xx})\mathrm{d}x, \quad
	H_2 = \int_R(u-u_{xx})v \mathrm{d}x,  \quad
	H_3 = \int_R u^3\mathrm{d}x,
\end{align*}
where $4v - v_{xx} = u$.
%These three types of equations share many common features such as complete integrability, Lax pairs and a bi-Hamiltonian structure \cite{ELY2006}.

Up to now, there have been thousands of papers focusing on the theoretical and numerical studies on these three equations.
It is worth mentioning that the invariant-preserving property is a key index of the success for numerical methods.
However, high-order invariants are usually difficult to preserve numerically. Liu \emph{et al.} also pointed out
``it appears a rather difficult task to preserve all three conservation laws'' in \cite{LX2016}. In this work, higher-order invariants of these equations will be re-derived in view of the energy method, which may be possible to provide some thoughts for invariant-preserving numerical methods. Actually, the energy method originated from conservation laws in physics was first proposed in 1928 by Courant, Friedrichs and Lewy \cite{CFL1928}. From then on, it has been widely applied to the mathematical and numerical analysis of nonlinear evolution equations. We trust the readers with \cite{Sun2018} instead of a long list of references to relevant works.

The rest of the paper is arranged as follows. In Section \ref{Sec_2}, combining the energy method and a skew-adjoint operator, we show the high-order invariants for the KdV equation, the Camassa--Holm equation and the Degasperis--Procesi equation, respectively. Then we list several applications for seeking some high-order invariants of other types of the shallow water wave equations in Section \ref{Sec_3}.

\section{Main results}\label{Sec_2}
In what follows, we directly show that $M_i$ $(i = 1,2,3,4)$, $E_i$ $(i = 1,2,3)$ and $H_i$ $(i = 1,2,3)$ are invariants of \eqref{eq1.2}, \eqref{eq1.3} and \eqref{eq1.4} subjected to the periodic boundary conditions based on the energy method, respectively.
\subsection{Invariants of the KdV equation}
\begin{proof}
(I) Multiplying by $1$, $u$ and $(u^2 + u_{xx})$, respectively, with \eqref{eq1.2}, we have $M_i$ $(i = 1,2,3)$.
In what follows, we show the fourth invariant $M_4$ of the KdV equation by the energy method.

Multiplying both sides of \eqref{eq1.2} by $2u_{xxxx} + \frac{10}{3}u_x^2+\frac{20}{3}uu_{xx}+\frac{20}{9}u^3$ and integrating the result, we have
\begin{align}
	%E'(t) &= \int_R \Big[2u_{xx}u_{xxt} -\frac{10}{3}(u_tu_{x}^{2}+2uu_xu_{xt})+\frac{20}{9}u^3u_t \Big]\mathrm{d}x \notag\\
%	&= \int_R\Big[2u_{xxxx}-\frac{10}{3}u_x^2 + \frac{20}{3}(uu_x)_x+\frac{20}{9}u^3\Big]u_t\mathrm{d}x  \notag\\
	0=& \int_R \Big(2u_{xxxx} + \frac{10}{3}u_x^2+\frac{20}{3}uu_{xx}+\frac{20}{9}u^3\Big)\cdot u_t\mathrm{d}x \notag\\
	  &- \int_R \Big(2u_{xxxx} + \frac{10}{3}u_x^2+\frac{20}{3}uu_{xx}+\frac{20}{9}u^3\Big)\cdot (4uu_x + 2u_{xxx})\mathrm{d}x  \notag\\
=& \int_R \Big[2u_{xx}u_{xxt} -\frac{10}{3}(u_tu_{x}^{2}+2uu_xu_{xt})+\frac{20}{9}u^3u_t \Big]\mathrm{d}x \notag\\
&- \int_R \Big(2u_{xxxx} + \frac{10}{3}u_x^2+\frac{20}{3}uu_{xx}+\frac{20}{9}u^3\Big)\cdot (4uu_x + 2u_{xxx})\mathrm{d}x  \notag\\
	=&\frac{{\rm d}}{{\rm d}t}M_4 -  8\int_Ruu_xu_{xxxx}\mathrm{d}x - \frac{40}{3}\int_R uu_{x}^3\mathrm{d}x - \frac{80}{3}\int_R u^2u_xu_{xx}\mathrm{d}x - \frac{80}{9}\int_Ru^4u_x\mathrm{d}x  \notag\\
	& - 4\int_Ru_{xxx}u_{xxxx}\mathrm{d}x - \frac{20}{3}\int_Ru_{xxx}u_{x}^2\mathrm{d}x - \frac{40}{3}\int_Ruu_{xx}u_{xxx}\mathrm{d}x - \frac{40}{9}\int_Ru^3u_{xxx}\mathrm{d}x.  \label{eq1.18}
\end{align}
It remains to check that the sum of all the integral terms in the above equation is zero. Calculating each term in \eqref{eq1.18} using the integration by parts, we have
\begin{align}
&-8\int_Ruu_xu_{xxxx}\mathrm{d}x = -20\int_Ru_xu_{xx}^2\mathrm{d}x,  \label{eq1.19}\\
&- \frac{80}{3}\int_R u^2u_xu_{xx}\mathrm{d}x =  \frac{80}{3}\int_Ruu_x^3\mathrm{d}x, \label{eq1.20}\\	
&- \frac{80}{9}\int_Ru^4u_x\mathrm{d}x = 0,\label{eq1.21}\\
&-4\int_Ru_{xxx}u_{xxxx}\mathrm{d}x = 0, \label{eq1.22}
\end{align}
\begin{align}
&- \frac{20}{3}\int_Ru_{xxx}u_{x}^2\mathrm{d}x =  \frac{40}{3}\int_Ru_xu_{xx}^2\mathrm{d}x,\label{eq1.23}\\
&- \frac{40}{3}\int_Ruu_{xx}u_{xxx}\mathrm{d}x = \frac{20}{3}\int_Ru_xu_{xx}^2\mathrm{d}x,\label{eq1.24}\\
&-\frac{40}{9}\int_Ru^3u_{xxx}\mathrm{d}x = -\frac{40}{3}\int_Ruu_x^3\mathrm{d}x. \label{eq1.25}
\end{align}
Substituting \eqref{eq1.19}--\eqref{eq1.25} into \eqref{eq1.18}, we have
	$\frac{{\rm d}}{{\rm d}t}M_4 = 0$,
which completes the proof.
\end{proof}
\begin{remark}
Suppose the general form of the {KdV equation} is
\begin{align*}
	u_t - auu_x - bu_{xxx} = 0,
\end{align*}
and the corresponding high-order invariant
\begin{align*}
M(t) = \int_R (u_{xx}^2 - Auu_x^2 + Bu^4)\mathrm{d}x.
\end{align*}
Using the same method above, we could derive
\begin{align*}
\left\{
\begin{array}{ll}
5a = 3Ab, \\
12Bb = Aa,
\end{array}
\right.
\end{align*}
which can be rewritten as
\begin{align*}
\frac{a}{b} = \frac{3A}{5} = \frac{12B}{A}.
\end{align*}
Therefore, it follows
\begin{align*}
A^2 = 20B.
\end{align*}
For instance, when $a=-6$, $b=-1$, we have
$A=10$, $B = 5$,
which deduces to the KdV equation as
\begin{align*}
	u_t + 6uu_x + u_{xxx} = 0,
\end{align*}
with a fourth-order invariant
\begin{align*}
M(t) = \int_R (u_{xx}^2 - 10 uu_x^2 + 5u^4) \mathrm{d}x.
\end{align*}
\end{remark}

\subsection{Invariants of the Camassa--Holm equation}
\begin{proof}
Multiplying by $1$ and $u$ on both sides of \eqref{eq1.3}, respectively, and then integrating the results, which implies $E_1$ and $E_2$  through the integration by parts. Below, we prove $E_3$ by the energy method. Firstly, noticing that \eqref{eq1.3} can be written with  a skew-adjoint operator $(1-\partial_{xx})^{-1}$ as
\begin{align*}
u_t + uu_x + \partial_x(1-\partial_{xx})^{-1}\Big(u^2+\frac12 u_x^2\Big) = 0.
\end{align*}
Let $g = (1-\partial_{xx})^{-1} \Big(u^2+\frac12 u_x^2\Big)$. Then we see from the above equation that \eqref{eq1.3} is equivalent to
\bcase
u_t + uu_x + g_x = 0,   \label{eq1.5a} \\
g - g_{xx} = u^2 + \frac{1}{2}u_x^2.  \label{eq1.5b}
\ecase
Multiplying \eqref{eq1.5a} by $3u^2 + u_x^2 - 2(uu_x)_x$ and integrating the result on both sides, we have
\begin{align}
0 &= \int_R(u_t + uu_x + g_x)\cdot ( 3u^2 + u_x^2 - 2(uu_x)_x )\mathrm{d}x  \notag\\
  &= \int_R u_t\cdot (3u^2 +  u_x^2 - 2(uu_x)_x )\mathrm{d}x + \int_R (uu_x+g_x)\cdot ( 3u^2 + u_x^2 - 2(uu_x)_x )\mathrm{d}x\notag\\
  &\triangleq A + B.  \label{eq1.7}
\end{align}
Calculating each term derives that
\begin{align}
A &= \int_R u_t\cdot (3u^2 +  u_x^2 - 2(uu_x)_x) \mathrm{d}x  \notag\\
&= \int_R u_t\cdot (3u^2 +  u_x^2)\mathrm{d}x + \int_R 2uu_x\cdot u_{xt}\mathrm{d}x  \notag\\
&= \int_R u_t \cdot 3u^2 \mathrm{d}x + \int_R u_t\cdot u_x^2 \mathrm{d}x + \int_R u \cdot (u_x^2)_t \mathrm{d}x \notag\\
&= \int_R (u^3)_t \mathrm{d}x + \int_R (u\cdot u_x^2)_t \mathrm{d}x \notag\\
&= \frac{\rm d}{{\rm d}t}\int_R(u^3 + uu_x^2) \mathrm{d}x    \label{eq1.8}
\end{align}
and
\begin{align}
B &= \int_R (uu_x+ g_x)\cdot ( 3u^2 + u_x^2 - 2(uu_x)_x ) \mathrm{d}x  \notag\\
&= \int_R u \cdot u_x^3 \mathrm{d}x + \int_R g_x \cdot (3u^2+u_x^2) \mathrm{d}x - \int_R g_x \cdot 2(uu_x)_x \mathrm{d}x  \notag\\
&= \int_R u \cdot u_x^3\mathrm{d}x + \int_R g_x \cdot( 3u^2+u_x^2) \mathrm{d}x + 2\int_R g_{xx} \cdot uu_x \mathrm{d}x  \notag\\
&= \int_R u\cdot u_x^3 \mathrm{d}x + \int_R g_x \cdot( 3u^2+u_x^2)\mathrm{d}x + 2\int_R(g-u^2-\frac{1}{2}u_x^2)\cdot uu_x \mathrm{d}x  \notag\\
&= \int_R g_x \cdot( 3u^2+u_x^2) \mathrm{d}x + 2\int_R g \cdot uu_x \mathrm{d}x  \notag\\
%&= \int_R g_x \cdot( 3u^2+u_x^2) \mathrm{d}x + \int_R g \cdot (u^2)_x \mathrm{d}x  \notag\\
&= \int_R g_x \cdot( 3u^2+u_x^2) \mathrm{d}x - \int_R g_x \cdot u^2 \mathrm{d}x   \notag\\
&= \int_R g_x \cdot( 2u^2+u_x^2) \mathrm{d}x \notag\\
&= 2 \int_R g_x \cdot (g-g_{xx}) \mathrm{d}x =0.    \label{eq1.9}
\end{align}
Substituting \eqref{eq1.8} and \eqref{eq1.9} into \eqref{eq1.7}, we have
\begin{align*}
\frac{\rm d}{{\rm d}t} \int_R (u^3+uu_x^2) \mathrm{d}x = 0,
\end{align*}
which implies $E_3$.
%\begin{align*}
%\int_R (u^3 + uu_x^2) \mathrm{d}x  \equiv Const.
%\end{align*}
\end{proof}

\subsection{Invariants of the Degasperis--Procesi equation}
\begin{proof}
Integrating on both sides of \eqref{eq1.4}, it easily obtains $H_1$. Then we show invariants $H_2$ and $H_3$ of \eqref{eq1.4}, respectively.
Firstly let $g = (1-\partial_{xx})^{-1} \Big(\frac32 u^2\Big)$, then \eqref{eq1.4} is equivalent to
\bcase
u_t + uu_x + g_x = 0,  \label{eq1.10} \\
g - g_{xx} = \frac{3}{2} u^2.  \label{eq1.11}
\ecase
Multiplying by $2u-6v$ on both sides of \eqref{eq1.10} and then integrating the result, we have
\begin{align}
0 &= \int_R (u_t + uu_x + g_x)\cdot ( 2u-6v)  \mathrm{d}x    \notag \\
  &= \int_R u_t \cdot(2u-6v) \mathrm{d}x + \int_R uu_x \cdot(2u-6v)  \mathrm{d}x + \int_R g_x\cdot(2u-6v) \mathrm{d}x  \notag \\
  &\triangleq C + D.   \label{eq1.12}
\end{align}
The each term in the above identity is estimated as
\begin{align}
C &= \int_R u_t \cdot(2u-6v)  \mathrm{d}x  = 2\int_R u_t \cdot u  \mathrm{d}x - 6\int_R u_t \cdot v \mathrm{d}x = 2\int_R u_t \cdot u \mathrm{d}x - 6\int_R (4v_t - v_{xxt}) \cdot v  \mathrm{d}x  \notag\\
  &= 2\int_R u_t \cdot u \mathrm{d}x - 24\int_R v_t \cdot v \mathrm{d}x - 6\int_R v_{xt} \cdot v_x  \mathrm{d}x   = \frac{{\rm d}}{{\rm d}t} \int_R ( u^2 - 12v^2 - 3v_x^2  ) \mathrm{d}x \notag\\
  &= \frac{{\rm d}}{{\rm d}t} \int_R\Big( u^2 - 3(4v-v_{xx})\cdot v  \Big)  \mathrm{d}x  = \frac{{\rm d}}{{\rm d}t} \int_R ( u^2 - 3uv  )  \mathrm{d}x   = \frac{{\rm d}}{{\rm d}t} \int_R u\cdot(u-3v)  \mathrm{d}x  \notag\\
  &= \frac{{\rm d}}{{\rm d}t} \int_R u\cdot(v-v_{xx}) \mathrm{d}x   = \frac{{\rm d}}{{\rm d}t} \int_R (u-u_{xx})\cdot  v \mathrm{d}x   \label{eq1.13}
\end{align}
and
\begin{align}
D &= \int_R uu_x\cdot(2u-6v)\mathrm{d}x +\int_R g_x \cdot (2u-6v)\mathrm{d}x  \notag\\
&= -6\int_R uu_x\cdot v \mathrm{d}x + \int_R g_x\cdot (2u-6v)\mathrm{d}x \notag\\
%&= -3\int_R (u^2)_x\cdot v\mathrm{d}x  + \int_R g_x\cdot (2u-6v)\mathrm{d}x \notag\\
&= 3\int_R u^2\cdot v_x \mathrm{d}x  + \int_R g_x\cdot (2u-6v) \mathrm{d}x \notag\\
&= 2\int_R (g-g_{xx})\cdot v_x \mathrm{d}x + \int_R g_x\cdot(2u-6v)\mathrm{d}x \notag\\
&= 2\int_R g \cdot v_x \mathrm{d}x - 2\int_R g_{xx}\cdot v_x\mathrm{d}x + \int_R g_x \cdot (2v-2v_{xx})\mathrm{d}x \notag\\
&= 2\int_R g\cdot v_x \mathrm{d}x + 2\int_R g_x \cdot v\mathrm{d}x - 2\int_R g_{xx}\cdot v_x \mathrm{d}x - 2\int_R g_x \cdot v_{xx} \mathrm{d}x \notag\\
&= 2\int_R (g v)_x \mathrm{d}x - 2\int_R(g_x \cdot v_x)_x\mathrm{d}x =0 .  \label{eq1.14}
\end{align}
Substituting \eqref{eq1.13} and \eqref{eq1.14} into \eqref{eq1.12}, we have
\begin{align*}
\frac{{\rm d}}{{\rm d}t} \int_R (u-u_{xx})\cdot v \mathrm{d}x = 0,
\end{align*}
which implies $H_2$.
%\begin{align*}
%\int_R (u-u_{xx})v\mathrm{d}x \equiv Const.
%\end{align*}

Finally, we show $H_3$. %Noticing \eqref{eq1.4} and \eqref{eq1.10}, we have
%\begin{align*}
%	g_x = -u_{xxt} + 3uu_x - (uu_x)_{xx}.
%\end{align*}
Multiplying \eqref{eq1.10} on both sides by $u^2$ and integrating the result, it yields by noting \eqref{eq1.11}
\begin{align*}
0 &= \int_R(u_t + uu_x + g_x) \cdot u^2 \mathrm{d}x  \notag\\
  &= \int_R u_t \cdot u^2 \mathrm{d}x + \int_R u^3\cdot u_x \mathrm{d}x+ \int_R g_x \cdot u^2 \mathrm{d}x  \notag\\
  &= \int_R \Big(\frac13 u^3\Big)_t \mathrm{d}x + \frac23\int g_x\cdot (g-g_{xx}) \mathrm{d}x\notag\\
  &= \frac13  \frac{{\rm d}}{{\rm d}t} \int_R u^3  \mathrm{d}x,
%0 &= (u_t + uu_x + g_x, u^2)  \notag\\
%&= (u_t,u^2) + (uu_x,u^2) + (g_x,u^2)  \notag\\
%&= \frac{1}{3} (u^3,1)_t + (-u_{xxt} + 3uu_x - (uu_x)_{xx},u^2)  \notag\\
%&= \frac{1}{3} (u^3,1)_t - (u_t, (u^2)_{xx}) - (uu_x, (u^2)_{xx})  \notag\\
%&= \frac{1}{3} (u^3,1)_t - (u_t + uu_x, (u^2)_{xx})  \notag\\
%&= \frac{1}{3} (u^3,1)_t + (g_x, \frac{2}{3}(g-g_{xx})_{xx}) \notag\\
%&= \frac{1}{3} (u^3,1)_t + \frac{2}{3} (g_{xxx}, g-g_{xx})\notag\\
%&= \frac{1}{3} (u^3,1)_t.
\end{align*}
which implies the invariant $H_3$.
\section{Applications to other periodic nonlinear dispersive waves}\label{Sec_3}
\subsection{Benjamin-Bona-Mahony equation}
Consider the Benjamin-Bona-Mahony equation \cite{MM1977} of the form
\begin{align}
u_t - u_{xxt} + u_x + \varepsilon uu_x = 0,\quad x\in R.  \label{eq2.1}
\end{align}
It can be written as
\begin{align*}
u_t + \partial_x(1-\partial_{xx})^{-1}\Big(u + \frac{\varepsilon}{2}u^2\Big) = 0,\quad x\in R.
\end{align*}
Let $g = (1-\partial_{xx})^{-1}\Big(u+\frac{\varepsilon}{2}u^2\Big)$, then the equation \eqref{eq2.1} turns out to be
\bcase
u_t + g_x  = 0, \label{eq2.2} \\
g -g_{xx} = u + \frac{\varepsilon}{2}u^2.  \label{eq2.3}
\ecase
Multiplying both sides of \eqref{eq2.2} by $u^2$ and integrating the result, and then using \eqref{eq2.3}, we have
\begin{align*}
0 &= \int_R (u_t+g_x)\cdot u^2 \mathrm{d}x   = \int_R u_t\cdot u^2 \mathrm{d}x + \int_R g_x\cdot u^2 \mathrm{d}x \notag \\
&= \int_R u_t\cdot u^2 \mathrm{d}x  + \frac{2}{\varepsilon}\int_R g_x\cdot (g-g_{xx}-u) \mathrm{d}x  = \int_R u_t\cdot u^2 \mathrm{d}x -\frac{2}{\varepsilon}\int_R g_x \cdot u \mathrm{d}x  \notag\\
&= \int_R u_t \cdot u^2\mathrm{d}x + \frac{2}{\varepsilon}\int_R  u_t\cdot u \mathrm{d}x
%= \frac{1}{3}\int_R(u^3)_t\mathrm{d}x + \frac{1}{\varepsilon}\int_R(u^2)_t \mathrm{d}x \notag\\
= \frac{{\rm d}}{{\rm d}t} \int_R \Big (\frac{1}{3}u^3 + \frac{1}{\varepsilon}u^2 \Big) \mathrm{d}x,  \\
\end{align*}
which indicates
\begin{align*}
\int_R \frac{1}{3}\Big(u^3 + \frac{1}{\varepsilon}u^2\Big) \mathrm{d}x
\end{align*}
is a three-order invariant for \eqref{eq2.1}.
\subsection{Regularized long wave equation}
Consider the regularized long wave equation \cite{SF1984} of the form
\begin{align}
u_t - \mu u_{xxt} + u_x + u^pu_x = 0,  \label{eq2.4}
\end{align}
where $\mu>0$ is a positive constant. When $p=2$, it is called modified regularized long wave equation; when $p\geqslant 3$, it is called generalized regularized long wave equation.
Similar to the foregoing argument, \eqref{eq2.4} can be written as an equivalent form of
\bcase
u_t + g_x = 0, \label{eq2.5}   \\
g-\mu g_{xx} = u+\frac{1}{p+1}u^{p+1}. \label{eq2.6}
\ecase
Multiplying both sides of \eqref{eq2.5} by $u^{p+1}$, integrating the result, and then using \eqref{eq2.6}, we have
\begin{align*}
0 &= \int_R(u_t + g_x)\cdot u^{p+1}\mathrm{d}x = \int_R u_t\cdot u^{p+1}\mathrm{d}x + \int_R g_x\cdot u^{p+1}\mathrm{d}x  \\
&= \int_R u_t\cdot u^{p+1} \mathrm{d}x + ({p+1})\int_R g_x\cdot (g-\mu g_{xx}-u)\mathrm{d}x \\
&= \int_R u_t \cdot u^{p+1}\mathrm{d}x - ({p+1})\int_R g_x\cdot u\mathrm{d}x  \\
&= \int_R u_t\cdot u^{p+1} \mathrm{d}x + ({p+1})\int_R u_t\cdot u\mathrm{d}x   \\
%&= \frac{1}{p+2}\int_R(u^{p+2})_t\mathrm{d}x + \frac{{p+1}}{2}\int_R(u^2)_t\mathrm{d}x  \\
&= \frac{{\rm d}}{{\rm d}t}\int_R \Big(\frac{1}{p+2}u^{p+2} + \frac{{p+1}}{2}u^2\Big)\mathrm{d}x,
\end{align*}
which indicates
\begin{align*}
\int_R \Big( \frac{1}{p+2}u^{p+2} + \frac{{p+1}}{2}u^2\Big) \mathrm{d}x
\end{align*}
is a high-order invariant for \eqref{eq2.4}. This corrects an invariant $I_3$ in Example 4 appeared in \cite{GO2018} (pp. 492).

\subsection{Rosenau equation}
Consider the Rosenau equation \cite{Park1990}
\begin{align}
u_t + u_{xxxxt} + u_x +uu_x = 0, \label{eq2.7}
\end{align}
which is equivalent to
\bcase
u_t + g_x = 0, \label{eq2.8} \\
g + g_{xxxx} = u + \frac{1}{2}u^2. \label{eq2.9}
\ecase
Multiplying both sides of \eqref{eq2.8} by $u^2$ and noticing \eqref{eq2.9}, similar to the argument in the above, we have a third-order invariant for \eqref{eq2.7} of the form
%\begin{align*}
%0 &= (u_t + g_x, u^2) \\
%&= (u_t,u^2) + (g_x,u^2) \\
%&= (u_t,u^2) + 2(g_x,g+g_{xxxx}-u)  \\
%&= (u_t,u^2) - 2(g_x,u) \\
%&= (u_t,u^2) + 2(u_t,u) \\
%&= \frac{1}{3}(u^3,1)_t + (u^2,1)_t \\
%&= (\frac{1}{3}u^3 + u^2,1)_t,
%\end{align*}
%which indicates
\begin{align*}
\int_R \Big(\frac{1}{3}u^3 + u^2\Big) \mathrm{d}x.
\end{align*}
\end{proof}

\section*{Acknowledgement}
We appreciate Prof. Zhi-zhong Sun for  many useful discussions.
This work is dedicated to Prof. Zhi-zhong Sun on the occasion of his 60th birthday.
The work is supported by Natural Science Foundation of Zhejiang Province (Grant No. LZ23A010007).

\end{document}